\documentclass[conference]{IEEEtran}
\IEEEoverridecommandlockouts
\usepackage{cite}
\usepackage{amsmath,amssymb,amsfonts}
\usepackage{algorithmic}
\usepackage{graphicx}

\usepackage{float}
\usepackage{textcomp}
\usepackage{xcolor}
\def\BibTeX{{\rm B\kern-.05em{\sc i\kern-.025em b}\kern-.08em
    T\kern-.1667em\lower.7ex\hbox{E}\kern-.125emX}}
\begin{document}

\title{Stochastic Economic Dispatch Considering Demand Response and
Endogenous Uncertainty\\
}

\author{\IEEEauthorblockN{Nasrin Bayat, Qifeng Li}
\IEEEauthorblockA{\textit{Dept. of Electrical and Computer Engineering} \\
\textit{University of Central Florida}\\
Orlando, USA\\
nasrinbayat@knights.ucf.edu}
\and
\IEEEauthorblockN{Joon-Hyuk Park}
\IEEEauthorblockA{\textit{Dept. of Mechanical and Aerospace Engineering} \\
\textit{University of Central Florida}\\
Orlando, USA \\
JoonPark@ucf.edu}
}
\maketitle
\begin{abstract}
This paper considers endogenous uncertainty (EnU) in the stochastic economic dispatch (SED) problem, where the endogenous uncertainty means decision dependent uncertainty. In this problem, demand response (DR) commitment is the source of the EnU. Nevertheless, EnU is not well considered in existing literature. Our first contribution is to build up an optimization model of DR-involved SED under EnU (SED-DR-EnU). This is a computational challenging problem due to the EnU. Our second contribution is introducing a coupled learning enabled optimization algorithm which can effectively solve the proposed SED-DR-EnU problem. This strategy is tested on the IEEE 14 bus, and IEEE 39 bus systems, and the results showed the importance of considering EnU in the DR-involved SED problem.
\end{abstract}
\section{Introduction}
Economic dispatch (ED) is the process of allocating output (active and reactive power) among available generation at the lowest cost and meeting the operational restrictions of the generation and the transmission system \cite{b1}. Demand response (DR) is a tariff or program that aims to persuade end-user consumers to change their electricity consumption in response to changes in the cost of electricity over time, or to offer incentive payments designed to persuade consumers to use less electricity during times of high market prices or when the stability of the grid is in danger \cite{b2}.

In reality, several variables, including inconvenience costs, incentive policies, and level of education, may affect how responsive consumers are. We must consider how dependent the consumers' uncertain behavior is on these elements.
The assumption that consumers' demand responsiveness or demand response commitment, is a constant that can be comprehended or perfectly projected by the utility company in advance, is made in the majority of recent publications \cite{b3}, \cite{b4}, \cite{b5}, \cite{b6}.
consumers may have different consumption patterns and preferences, making it challenging for utility companies to understand each DR consumer's unique characteristics in real-world settings.
Additionally, even if such data were available, other factors, such as unusual events, might still affect people's responsiveness. A deterministic model would find it difficult to forecast the effects of these factors. As a result, while the system is in operation, the utility companies should be completely unaware of the system's actual DR capacity. However, while plans were being created, these DR ambiguity issues were rarely taken into consideration.

Exogenous uncertainty or decision-independent uncertainty is the term used in practice to describe an uncertain variable that follows a particular statistical regularity that may be known in advance and does not change over time. On the other hand, endogenous uncertainty (EnU) or decision-dependent uncertainty refers to uncertain variables whose probability distribution is not known and is affected by the decision variables. Neglecting EnU in DR modeling could lead utility companies to make poor system strategic planning decisions by vastly overestimating the capacity value of DR resources. 

In the prior works, the majority of the DR uncertainties were modeled as fixed probability distributions that could be fully described before the evaluation. In other words, existing models have considered consumer demand responsiveness as an exogenous uncertain quantity that is unrelated to the operator's choice of control method. In practice, consumer participation in DR might benefit them financially but adversely affect their comfort and accessibility of electricity service \cite{b7}. 
Ineffective scheduling strategies could reduce the level of compliance of the consumers with DR calls \cite{b8}. As a result, the current distribution of consumer responsiveness is influenced by grid operating decisions.
DR presents uncertainty and is dynamic in nature. Thus, underestimating these aspects is likely to yield unreliable utilization of DR-based economic dispatch (ED) \cite{b9}.
For example, the operational and control tactics may impact some stochastic features (such as response probability and temperature preference); however, these kinds of uncertainties, which are called EnU are disregarded and considered to be static and known probability distributions \cite{b10}.

 Stochastic programming (SP) methods use the discretization process with scenarios of uncertainty without knowing the precise probability distribution of the uncertainty. Scenarios are frequently regarded as a priori in SP approaches when the uncertainty is decision independent, which is fair because they may be simulated using either an exogenous prediction model or from a historical data collection \cite{b11}. 
Decisions taken before the realization of uncertainty affect EnU. By altering its parameters or information structures, decisions might have an impact on uncertainty. A few examples include operating budgets, equipment switch status, and investment plans for renewable energy sources (RES). End-consumers must adapt their loads when the demand response set point changes. Some examples include reducing PV output, switching the state of domestic appliances, and rescheduling electric vehicle charging. Decision-related uncertainty loads can change as a result of amplification or reduction of uncertainty sources \cite{b12}.

In this paper, the decision variable influences the probability distribution of the random variable. In addition, the learning step of predicting the latent dependency and the optimization step of computing a candidate decision are carried out interactively.
The contributions of this work are as follows:
\begin{itemize}
    \item Considering both the exogenous and endogenous uncertainties of Demand Response in the stochastic economic dispatch (SED-DR-EnU) problem, which are a result of demand-side physical factors and the behavior of consumers.
    \item A coupled learning enabled optimization (CLEO) algorithm \cite{b14} is introduced to effectively solve the DR problem under both exogenous and endogenous uncertainties without making any assumptions about the probability distribution of decision-dependent uncertainty.
    \item Comparative analysis to demonstrate how the demand-side participation uncertainty may have a significant impact on the effectiveness of DR programs. 
\end{itemize}

\section{Optimization Model of SED-DR-EnU}
This section describes the optimization model we developed for the stochastic economic dispatch considering demand response under both exogenous and endogenous uncertainties in detail.
Incorporating demand response commitments $P^{rd}_{j,max}$ from end-consumers utilizing a variety of incentives, a demand response provider (DRP) aggregates a group of residential/commercial consumers and sets its offer price $\pi^{DR}_j$ in the day-ahead market. 
The cost function of the unit i is presented as:
\begin{equation}
    C_i(P_{i}^{G}) = \sum_{i=1}^{N_G}\left(a_{i} {P_{i}^{G}}^2+b_{i} P^G_{i}\right)
\end{equation}

Where, $a_i$ and $b_i$ are coefficients of the cost function and $P_i^G$ is power generation of unit i. Each RES output's uncertainty interval is given by $\zeta_{i} \in\left[-P_{R i}, P_{R i}\right] \times(r / 100)$, where $P_{R i}$ is the $i_{th}$ RES nominal output and $r$ determines the uncertainty level of $P_{R i}$. It is assumed that the predicted RES generation is such that the average value of its uncertainty is zero.
The parameter $r$ shows the confidence interval for the RES output forecast error; that means, the closer $r$ is to 100, the power forecast would be the less accurate.
The power generation is made up of the nominal conditions output $P^G_{i(\text { base })}$, and the shift in generation brought on by RES output changes.
Let $\alpha_i$ stand for the unit $i$'s participation factor, which represents the contribution of generator $i$ to the uncertain component of total load, i.e., the difference between the real and forecasted total load. It is subject to the constraints (2) \cite{b15}.
\begin{equation}
\resizebox{.99\hsize}{!}{$
\alpha_{i}=\frac{1 / a_{ i}}{\sum_{j=1}^{N_G} 1 / a_{ j}},
\quad 0\leq \alpha_{i} \leq 1, \quad 
\sum_{j=1}^{N_G} \alpha_i = 1, \quad i=1, \ldots, N_G$}
\end{equation}
Then $P^G_i$ is calculated as:
\begin{equation}
P^G_{i}=P^G_{i(\text { base })}-\alpha_{i} \times \sum_{j=1}^{N_{R}} \zeta_{j}, \quad i=1, \ldots, N_{G}
\end{equation}

$N_G$ is the total number of generators and $N_R$ is the number of available RES. If $\lambda$ is covariance matrix of the predicted RES generation uncertainty and the expected prediction error of the uncertainty is zero,
the mathematical formulation of the expectation of supply-side cost, which is the sum of the generation costs, and the cost of DRPs will be as follows:
\begin{equation}
\begin{aligned}
\min \quad f\left( P^{r d}\right)= \mathbb{E}[C_i(P_{i}^{G}) + \sum_{i=1}^{N_G} a_{i}^{\prime} \alpha_{i}^{2} +\sum_{j=1}^{N_{DRP}}  P_{j,EnU}^{r d} \pi_j^{DR}] 
\end{aligned}
\end{equation}
Where 
\begin{equation}
a_{i}^{\prime}=\left[\sum_{j=1}^{N_{R}} \sum_{k=1}^{N_{R}} \lambda_{(j, k)}\right] \times a_{i}, i=1, \ldots, N_G
\end{equation}
\begin{equation}
   P_{j,EnU}^{r d}=h(P_{j}^{r d},\varepsilon_j)
\end{equation}
where $P_{j}^{r d}$ is the decision variable which represents the accepted DR commitment of the j-th DRP, $\varepsilon_j$ is a random variable that is independent of $P_{j}^{r d}$. $P_{j, EnU}^{r d}$ is the actual DR commitment considering EnU that shows the responsiveness of consumers and is a function of the independent system operator's (ISO) decisions (accepted DR commitments).  $\pi^{DR}_{j}$ is the j-th DRP's offer price, and $N_{DRP}$ is the number of DRPs for a group of consumers.
Subject to the following constraints:\\
1- Energy adequacy constraint :

\begin{equation}
\resizebox{.99\hsize}{!}{$
\sum_{i=1}^{N_G} P^G_{i}+\sum_{i=1}^{N_R} P^R_{i}+\sum_{j=1}^{N_{DRP}} P_{j,EnU}^{r d}>\sum_{i=1}^{N_B} P^L_{i}$}
\end{equation}
2- Power flow limit:
\begin{equation}
\resizebox{.99\hsize}{!}{$
L \times \left( I_{1}\sum_{i=1}^{N_G}P^G_{i}+I_{2}\sum_{j=1}^{N_{DRP}}   P_{j,EnU}^{r d}+I_{3}\sum_{i=1}^{N_R}\left( P^R_{i}+\zeta_i\right)-I_{4}\sum_{i=1}^{N_B} P^L_{i}\right) \leq F_{l}^{\max }$}
\end{equation}
where $L$ is an $N_{L} \times\left(N_{B}-1\right)$ matrix of power transfer distribution factors. $N_B$ is the number of buses and $N_L$ is the number of lines. $P^L_{i}$ represents the load at  bus $i$. $\epsilon$ indicates the upper
limit of violation probability in (7). $F_{l}^{\max }$ is the maximum line power flow vector. $I_{1}, I_{2}, I_{3}$, and $I_{4}$ are matrices that respectively map the vectors $P^{G},P^{rd}, P^{R}$, $\zeta$, and $P^{L}$ into $\left(N_{B}-1\right) \times 1$ vectors, with the nonzero elements denoting connections to one of the $N_{B}$ buses other than the slack bus.\\
3- Power generation limit:
\begin{equation}
P^G_{ i,min} \leq P^G_{i} \leq P^G_{ i,max}, i=1,2, \ldots, N_G
\end{equation}
4-  DR commitment limit:
\begin{equation}
0 \leq P_{\mathrm{j}}^{\mathrm{rd}} \leq P_{ \mathrm{j}, max}^{\mathrm{rd}}, j=1,2, \ldots, N_{DRP}
\end{equation}
 Assuming a linear aggregated demand curve, $P_{ \mathrm{j}, max}^{\mathrm{rd}}$ can be calculated as
 \begin{equation}
 P_{ \mathrm{j}, max}^{\mathrm{rd}}=\min \left(P_{\text {base }, j}, \frac{\pi_{s, j}}{\pi_{j}^{\max }-\pi^{R R}} P_{\text {base }, j}\right)
\end{equation}
where $\pi_{s,j} $ is DRP’s incentive price to end-consumers, $P_{base}$ is the baseline for end-consumers, $\pi_j^{max}$ is the demand curve's Y-axis intercept. $\pi^{RR}$ is the retail price.
This problem is very hard to solve, because $P_{j,EnU}^{r d}$ and the distribution of the random variable $\varepsilon$ are unknowns before optimization. In addition because realizations of $\varepsilon$ cannot be observed, traditional SP techniques cannot be used. Data pairs of decisions and uncertainties should be gathered together to understand the link between uncertainty and the decision variable in order to successfully address such SP problems.
This paper introduces a coupled learning enabled optimization (CLEO) algorithm \cite{b14} to solve this problem.
\section{Solution Method}
This section is devoted to providing a detailed description of the presented framework. 
In this work, it is considered that the responsiveness of consumers is a function of the independent system operator's (ISO) decisions as below: 
 \begin{equation}
     P_{j,EnU}^{r d} = \psi ( P_{j}^{r d}) + \varepsilon
 \end{equation}
 \begin{equation}
     \psi( P_{j}^{r d})= \left(A^{1}\right)^{\top} P_{j}^{r d}+\left(A^{0}\right)^{\top}
 \end{equation}
The derivative-free trust region method is the foundation of the CLEO algorithm, which is used to solve SED with latently decision-dependent uncertainty. In the proposed solution process, the SED-DR-EnU problem (1)-(9) is solved in two steps, learning and optimization. At each iteration, a region around the current best solution is defined that is trusted to be a sufficient approximation of the objective function and hence the model is called trust region. Using local linear regression (LLR) centered at the current iteration, the local latent dependency is predicted in the learning step. During the optimization phase, a potential answer to the trust region subproblem is sought, which is made up of the random LLR model. At the $k_{th}$ iteration, the current point is $P_{j,k}^{r d}$, $\delta_{k}$ is the radius of the current trust region and the data set is $T_{k}=\left\{\left(P_{j,i}^{r d}, P_{j,i,EnU}^{r d}\right)\right\}$ of size $N_k$.
 The estimation of parameters $A^1$ and $A^0$ can be found by  minimizing the sum of square errors between the actual and acctepted DR commitment.
\begin{equation}
\begin{aligned}
\left\{\widehat{A}^{k, 0}, \widehat{A}^{k, 1}\right\} \in \underset{A^{0}, A^{1}}{\operatorname{argmin}} \sum_{\left(P_{j,i}^{r d}, P_{j,i,EnU}^{r d}\right) \in T_{k}} \\ \left\|P_{j,i,EnU}^{r d}-\left(A^{1}\right)^{\top} P_{j,i}^{r d}-\left(A^{0}\right)^{\top}\right\|^{2}
\end{aligned}
\end{equation}
 $\widehat{A}^{k, 0}$ and $\widehat{A}^{k, 1}$ are the estimation of parameters $A^0$ and $A^1$ respectively. By subtracting the actual DR commitment from the estimated one, the error $e^{k, i}$ can be calculated as shown in \eqref{eq4}.
\begin{equation}
e^{k, i} \leftarrow P_{j,i,EnU}^{r d}-\left(\widehat{A}^{k, 1}\right)^{\top} P_{j,i}^{r d}-\left(\widehat{A}^{k, 0}\right)^{\top} \quad i=1,2,..,N_k
\label{eq4}
\end{equation}
A random variable with the empirical probability distribution of $\left\{e^{k, i}\right\}_{i=1}^{N_{k}}$ is denoted by the symbol $\varepsilon^{k}$. The $\mathrm{LLR}$ model $m_{k}$ is constructed as follows.
\begin{equation}
m_{k}\left(P_{k}^{r d}+s^k, \varepsilon^{k}\right) \triangleq\left(\widehat{A}^{k, 1}\right)^{\top}\left(P_{k}^{r d}+s^k\right)+\left(\widehat{A}^{k, 0}\right)^{\top}+\varepsilon^{k}
\end{equation}
Where $s^k$ is the step which minimizes the objective function in the trust region and satisfies a "sufficient" decrease requirement. The mathematical formulation of the DR problem under both exogenous and endogenous uncertainties (1)-(8) will change to: \\
\begin{equation}
\begin{aligned}
\min \quad f\left( P_{k}^{r d}+s^k\right)= \mathbb{E}[C_{k,i}(P_{k,i}^{G}) +\\ \sum_{i=1}^{N_G} a_{k,i}^{\prime} \alpha_{k,i}^{2} +\sum_{j=1}^{N_{DRP}}  m_{j,k}\left(P_{j,k}^{r d}+s^k, \varepsilon^{j,k}\right) \pi^{DR}_{j, k}] 
\end{aligned}
\end{equation}
subject to the following constraints:\\
1- Energy adequacy constraint: We reformulate the constraint (7) to (18) by substituting $P_{j,EnU}^{r d}$ with $m_{j,k} \left( P_{j,k}^{r d}+s^k,\varepsilon^{j,k}\right)$.
\begin{equation}
\resizebox{.99\hsize}{!}{$
\sum_{i=1}^{N_G} P^G_{k,i}+\sum_{i=1}^{N_R} P^R_{k,i}+\sum_{j=1}^{N_{DRP}} m_{j,k} \left( P_{j,k}^{r d}+s^k , \varepsilon^{j,k}\right)>\sum_{i=1}^{N_B} P^L_{i,k}$}
\end{equation}
2- Power flow limit:
\begin{equation}
\resizebox{.99\hsize}{!}{$
L \times \left( I_{1}\sum_{i=1}^{N_G}P^G_{k, i}+I_{2}\sum_{j=1}^{N_{DRP}}   P^{rd}_{ k,j}+I_{3}\sum_{i=1}^{N_R}\left( P^R_{k,i}+\zeta_i\right)-I_{4}\sum_{i=1}^{N_B} P^L_{i,k}\right) \leq F_{l}^{\max }$}
\end{equation}
3- Power generation limit:
\begin{equation}
P^G_{ i,min} \leq P^G_{k, i} \leq P^G_{ i,max}, i=1,2, \ldots, N_G
\end{equation}
4-  DR commitmen limit:
\begin{equation}
0 \leq P_{\mathrm{k},\mathrm{j}}^{\mathrm{rd}} \leq P_{ \mathrm{j}, max}^{\mathrm{rd}}, j=1,2, \ldots, N_{DRP}
\end{equation}
5- \begin{equation}
    \lVert s \rVert \leq \delta_k 
\end{equation}
Actual to predicted reduction ratio of the objective function is used to determine whether $s_k$ is a descent step. The problem has a decision-dependent uncertainty. This dependency is unknown. As a result, it is impossible to compute the objective function's true value. 
In this case, local linear regression models can be used to estimate the value of the objective function. Particularly, $u_{k}$ and $u_{k+1 / 2}$ are used to estimate the values of the functions $f\left(p^{k}\right)$ and $f\left(p^{k}+s^{k}\right)$.
\begin{equation}
\begin{aligned}
u_{k} & \triangleq \left[f\left(p^{k}, m_{k, 1}\left(p^{k}, \varepsilon^{k, 1}\right)\right)\right] 
\end{aligned}
\end{equation}
\begin{equation}
\begin{aligned}
u_{k+1 / 2} & \triangleq \left[f\left(p^{k}+s^{k}, m_{k, 2}\left(p^{k}+s^{k}, \varepsilon^{k, 2}\right)\right)\right]
\end{aligned}
\end{equation}
The estimated ratio between actual and predicted decrease in objective function is then roughly represented as
\begin{equation}
\rho_{k} \triangleq\left(u_{k}-u_{k+1 / 2}\right) /\left(f_{k}\left(p^{k}\right)-f_{k}\left(p^{k}+s^{k}\right)\right)
\end{equation}
To determine if the new step $p^{k}+s^{k}$ will be approved or denied, the ratio $\rho_{k}$ needs to be large enough. In addition, $\lVert s \rVert \leq \delta_k $ should be satisfied and the generalized gradient and the TR radius $\delta_{k}$ should be large enough in relation to one another.
After applying the two steps of the algorithm, cost and total demand response commitment are predicted.
\section{Results and Discussion}
In this section, IEEE 14 bus and 39 bus power systems with DRPs and RESs are examined to evaluate the performance of the proposed stochastic economic dispatch model. The SED-DR-ENU formulations are coded in Python.  Looking back on (11), a linear aggregated inherent demand curve for end-consumers is assumed, with baseline $Pbase = 60$, retail price $\pi^{RR}=100$, the j-th DRP’s offer price $\pi^{DR}_j = 100$, and Y-axis intercept $\pi^{max}_j=400$. In total, 3 cases are considered to analyze the effect of EnU on DR commitment, as it is described below:\\
Case 1: Demand response considering both endogenous and exogenous uncertainties (uncertainty of load and RES)\\
Case 2: Demand response without considering uncertainties. \\
Case 3: Demand response considering just exogenous uncertainties (uncertainty of RES)
\subsection{IEEE 14 bus system}
\begin{figure}[H]
    \centering
    \includegraphics[scale=0.5]{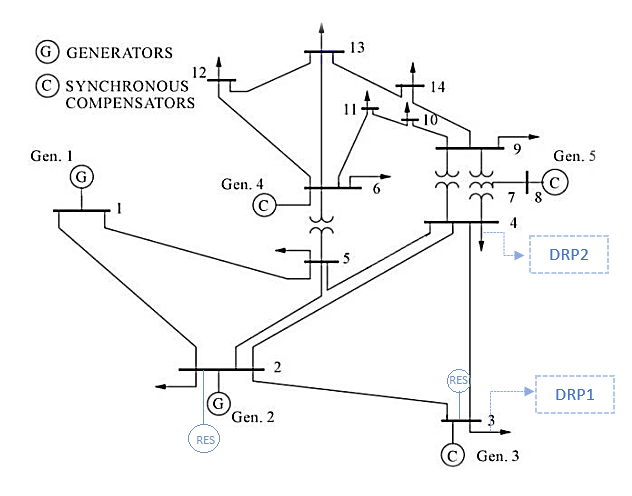}
    \caption{IEEE 14 Bus system with 2 DRPs and 2 RESs}
    \label{fig:14drp}
\end{figure}

\begin{figure}[h]
    \centering
    \includegraphics[scale=0.35]{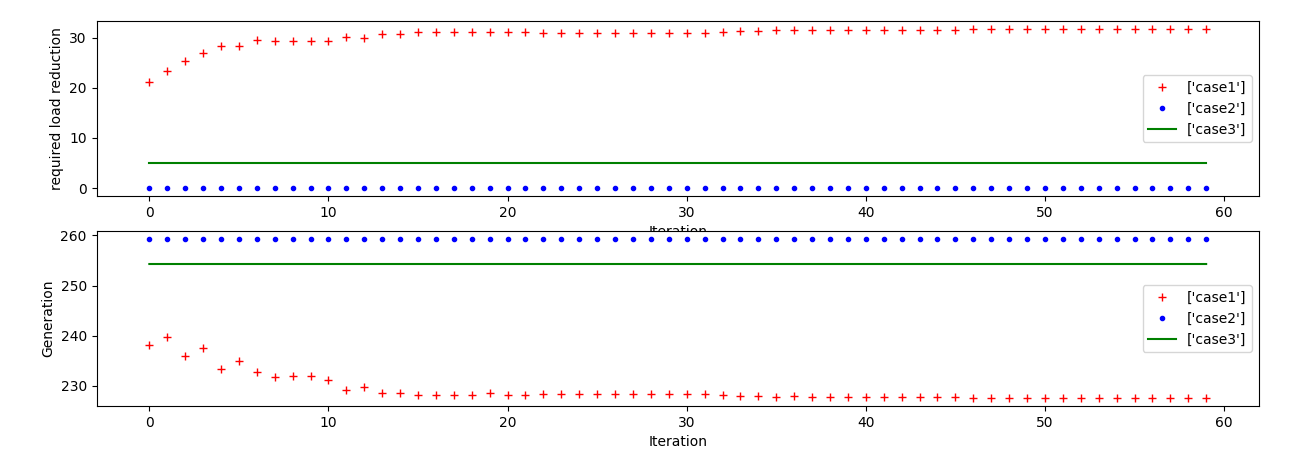}
    \caption{IEEE 14 Bus system results}
    \label{fig:14result}
\end{figure}

This subsection analyzes a more complex IEEE 14-bus system that has two DRPs and two RESs as shown in Fig.\ref{fig:14drp}. In order to lower the demand, it is supposed that DRPs are located in buses 3 and 4 that have the highest loads $(N_{DRP} = 2)$.
It is assumed that two RESs, each with a 40 MW capacity, are connected to buses 2 and 3, with $r = 20\%$ $(N_{R} = 2)$. Fig. \ref{fig:14result} shows the result of the simulation for the three cases. The amount of required load reduction for case 1 is higher than for cases 2 and 3 because case 1 is more accurate. Based on the first row of Table \ref{Tab:39143case} case 1 needs more DR commitment compared with cases 2 and 3. This result represents the importance of considering EnU in the SED problem. 
\subsection{IEEE 39 bus system}
\begin{figure}[H]
    \centering
    \includegraphics[scale=0.4]{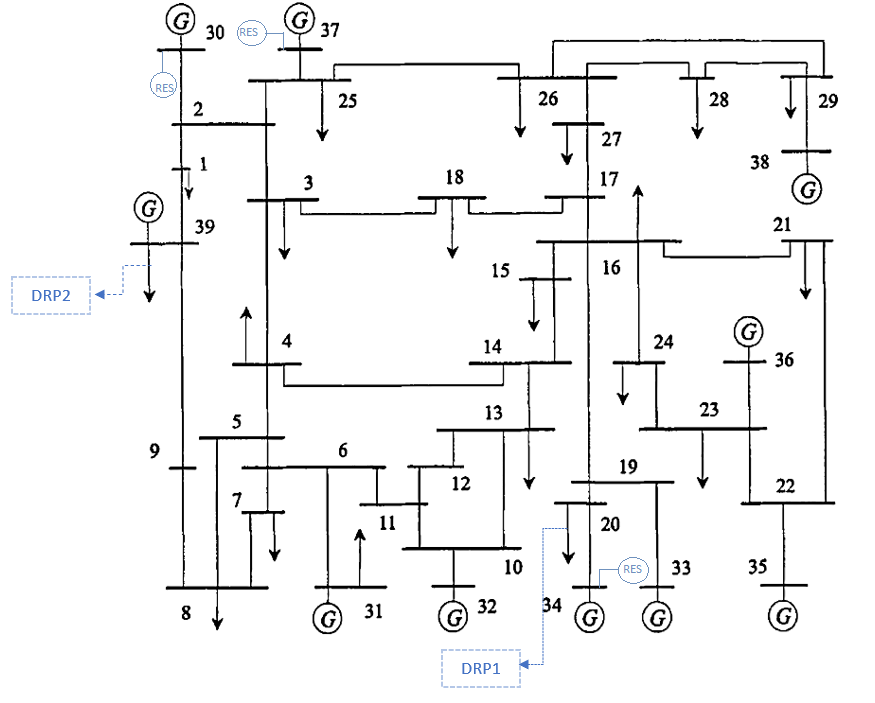}
    \caption{IEEE 39 Bus system with 2 DRPs and 2 RESs}
    \label{fig:39drp}
\end{figure}

\begin{figure}[h]
    \centering
    \includegraphics[width=\linewidth]{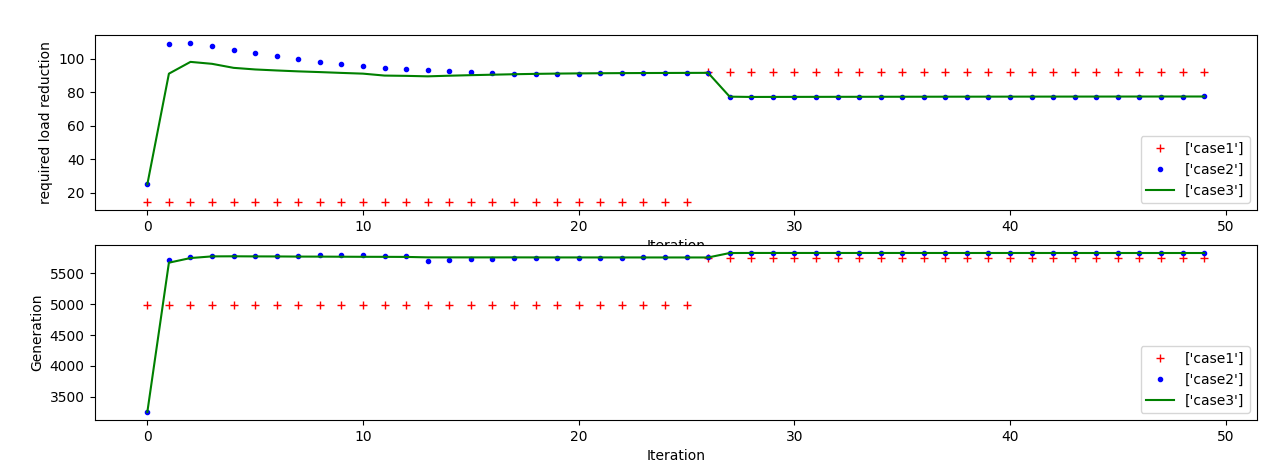}
    \caption{IEEE 39 Bus system results}
    \label{fig:39result}
\end{figure}

\begin{figure}[!h]
    \centering
    \includegraphics[scale=0.4]{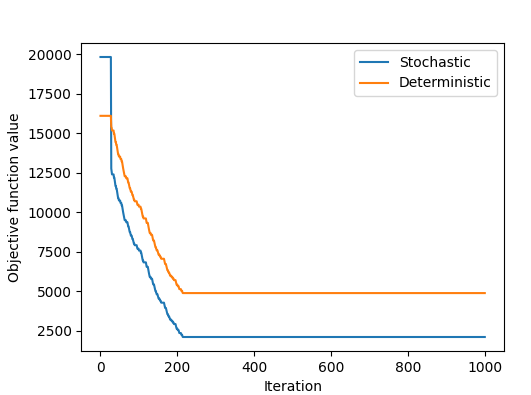}
    \includegraphics[scale=0.7]{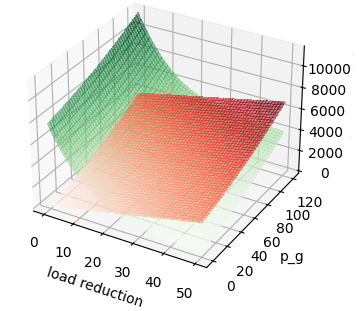}
    \caption{Comparison of convergence curves of stochastic (blue) and deterministic (orange) ED (top). Deterministic (red plot) VS. Stochastic (green plot) results. Vertical axis shows objective function vlue which has minimum value of 2109.73 and 4890 for stochastic and deterministic algorithms respectively. X, Y and Z axes show accepted DR
commitment, total power generation of the generators and objective function value respectively (bottom).}
    \label{fig:convcurves}
\end{figure}

In this subsection, IEEE 39 bus system with two DRPs is studied as shown in Fig. \ref{fig:39drp}. Buses 20 and 39 are presumed to have DRPs $(N_{DRP} = 2)$. RES integration is simulated using three wind farms connected to buses 30, 34 and 37. Total contribution of wind farms is assumed to be around 20$\%$ of the system's total generating capacity, $r$ is equal to 20$\%$ and $N_{R} = 3$. Fig. \ref{fig:39result} represents the result of the simulation for the three
cases. Required load reduction for case 1 is higher than cases 2 and 3, this shows the inaccuracy of the case 2 and 3. The second row of Table \ref{Tab:39143case} compares the DR commitment values for the 39 bus system.
Cases 2 and 3 represent the existing research where case 2 does not consider any uncertainties and case 3 considers just exogenous uncertainties (uncertainty of RES). Case 1 is the proposed scenario of this paper which is closer to real-world situations since it considers both endogenous and exogenous uncertainties (uncertainty of load and RES). Given this, the simulation results shown in \ref{fig:14result} and \ref{fig:39result} indicate that cases 2 and 3 (i.e., existing research) are not accurate enough to represent the real-world. An important conclusion we can make from the simulation results is that endogenous uncertainty should be considered in the economic dispatch problems which involve demand response.
To evaluate the performance of the stochastic ED model, the convergence curve of a deterministic ED which does not take into account any uncertainties is compared with our model in Fig. \ref{fig:convcurves}. Both approaches converge after about 210 iterations.
In addition, a comparison between stochastic and deterministic approaches is represented in Fig. \ref{fig:convcurves}. Minimum
objective function vlues are 2109.73 and 4890 for stochastic and deterministic
algorithms respectively.
The SED problem is also solved with SVR and linear regression in addition to CLEO algorithm, to evaluate CLEO performance.
Table \ref{Tab:svrcleo} compares the standard deviation of the objective function and simulation time of all three algorithms. Based on these results CLEO performs as well as the SVR algorithm.


\begin{table}[ht]
\caption{DR commitment value for 3 cases and IEEE 14 bus and 39 bus systems }
\centering
\begin{tabular}{|c|c|c|c|}
\hline
       Value & Case 1 & Case 2 & Case 3\\
         \hline
         DR commitment(14Bus)& 28.05 pu &6.48e-08 pu &5.30e-08 pu\\
         \hline
         DR commitment (39Bus)&100.44 pu &100.22 pu &100.26 pu\\
         \hline
    \end{tabular}
\label{Tab:39143case}
\end{table}
\begin{table}[ht]
\caption{Comparison of standard deviation of the supply-side cost and the  simulation time between SVR, CLEO and Linear Regression models.}
\centering
 \begin{tabular}{|c|c|c|c|}
 \hline
       Value &SVR & CLEO & LinearRegression\\
         \hline
         Obj func standard deviation& 0 &1.13e-13&0.10\\
         \hline
         Simulation time(1000 samples) &0.004&0.02&0.007\\
         \hline
    \end{tabular}
\label{Tab:svrcleo}
\end{table}
\section{Conclusion}
SED-DR-EnU framework is presented in this paper, to account for the importance of considering endogenous uncertainty in the day-ahead market. Demand response is provided by DRPs using bottom-up aggregation, which has a significant degree of DR uncertainty because of end-consumer behavior. 
The proposed framework's ability to take into account both of the exogenous and endogenous uncertainties of DR in the SED problem—which are caused by physical factors and behavior of consumers on the demand side—represents its key innovation in comparison to past research.
CLEO model has been presented to take into account the unpredictability of consumer responsiveness, and it is considered that the responsiveness of consumers is a function of accepted DR commitments, to show the dependence of consumer behavior with DR operations.
The simulation results show that DR commitment is influenced by factors like consumer consumption patterns and the grid's DR operation decisions.
Additionally, a comparative study demonstrates that demand-side participation uncertainty may have a high impact on how effectively DR programs operate. Thus, in actuality, effective calculation of DR commitment should take into account the impact of endogenous uncertainty.

\end{document}